\documentclass{amsart}
\usepackage{amsmath}
\usepackage{graphicx}
\usepackage[utf8]{inputenc}

\newtheorem{thm}{Theorem}[section]
\newtheorem*{thm*}{Theorem}

\newtheorem{lem}[thm]{Lemma}
\newtheorem{pro}[thm]{Proposition}
\theoremstyle{definition}

\theoremstyle{remark}
\newtheorem{rem}[thm]{Remark}

\numberwithin{equation}{section}
\newcommand{\norm}[1]{\left\Vert#1\right\Vert}
\newcommand{\abs}[1]{\left\vert#1\right\vert}
\newcommand{\set}[1]{\left\{#1\right\}}
\newcommand{\paren}[1]{\left(#1\right)}

\newcommand{\inv}[1]{\frac{1}{#1}}

\newcommand{\ca}[1]{\mathcal{#1}}

\newcommand{\ppp}{\ldots}

\newcommand{\Z}{\mathbb Z}

\newcommand{\p}{\mathcal P}
\newcommand{\h}{\mathcal H}
\newcommand{\m}{\mathcal M}

\newcommand{\e}{\varepsilon}
\newcommand{\la}{\lambda}

\begin{document}

\title[A Marstrand's projection type theorem]{A generalization of Marstrand's theorem}%
\author{Jorge Erick L\'{o}pez, Carlos Gustavo Moreira, Waliston Luiz Silva}%
\address{IMPA - Estrada D. Castorina, 110 - 22460-320 - Rio de Janeiro - RJ - Brasil}%
\email{jorgeerickl@gmail.com; gugu@impa.br; waliston@yahoo.com}%

\thanks{The first author was supported by the Balzan Research Project of J.Palis and by INCTMat and PCI projects. The second author was supported by the Balzan Research Project of J.Palis and by CNPq.}

\begin{abstract}
In this paper we prove two general results related to Marstrand's projection theorem in a quite general formulation over separable metric spaces under a suitable transversality hypothesis (the "projections" are in principle only measurable) - the result is flexible enough to, in particular, recover most of the classical Marstrand-like theorems. Our proofs use simpler tools than many classical works in the subject, where some techniques from harmonic analysis or special geometrical structures on the spaces are used.
\end{abstract}

\maketitle
\section{Introduction}
Let $X,Y$ be separable metric spaces, $(\Lambda, \p)$ a probability space and $\pi:\Lambda\times X\to Y$ a measurable function. Informally, one can think of $\pi_\la(\cdot) = \pi(\la,\cdot)$ as a family of projections parameterized by $\la$. We assume that for some positives real numbers $\alpha, \kappa$ and $C$ the following transversality property is satisfies:
\begin{equation}\label{transversality}
\p[\la\in \Lambda: d(\pi_{\la}(x_1),\pi_{\la}(x_2))\leq\delta d(x_1,x_2)^{\alpha}]\leq C\delta^{\kappa}
\end{equation}
for all $\delta>0$ and all $x_1,x_2\in X$. In most examples $\kappa=\dim Y$. 

We called the measure $\eta$ on $Y$ $\kappa$-\emph{regular} if it is Borel, upper regular by open sets and $$\liminf_{r\to0}\eta(B(y,r))/r^\kappa>0 \text{ for all }y\in Y.$$ Lebesgue measure in Euclidean spaces are regular strongly, in the sense that the limit is uniformly positive. In those cases we say that $\eta$ is \emph{strongly $\kappa$-regular}.

We are able to announce general versions of Marstrand theorem in this context. Assuming $X$ be an analytic subset of a complete separable metric space, then:
\begin{thm} \label{thm_Marstrand1} $\dim \pi_\la(X)\geq\min(\kappa,\dim X/\alpha)$ for a.e. $\la\in\Lambda$.
\end{thm}
In particular, if $\kappa=\dim Y$ and $\pi_\la$ are $\alpha$-H\"older when $\dim X<\alpha\kappa$, then $\dim\pi_\la(X)=\min(\kappa,\dim X/\alpha)$ for a.e. $\la\in\Lambda$.

\begin{thm}\label{thm_Marstrand2}
Suppose $\dim X>\alpha\kappa$.
\begin{itemize}
\item[(i)] If $\eta$ is $\kappa$-regular, then $\eta(\pi_\la(X))>0$ for a.e. $\la\in\Lambda$.
\item[(ii)] If $\eta$ is strongly $\kappa$-regular and $\sigma$-finite, then $\int_{\Lambda}\eta(\pi_\la(X))^{-1}d\p<+\infty$.
\end{itemize}
\end{thm}

\begin{rem}
Those Theorems recover the Solomyak formulation in \cite{Solomyak_1998} Theorem 5.1. Solomyak scheme has the interesting idea to consider the exponent $\alpha$ depending continuously on the $\lambda$-parameter. As Solomyak, we also use the ideas of Mattila in \cite{Mattila_1995} to prove the Marstrand's Theorem. However, some of them can not be extended to this general context. In those cases our modifications stay elementary enough.
\end{rem}
\section{Preliminaries}
For $X$ separable metric space, let

\begin{align*}
\ca{M}^s(X)=\{\mu: \mu \textrm{ Borel measure on $X$ with }0<\mu(X)<\infty \textrm{ such that there is $c>0$}&\\\text{such that }\mu(B(x,r))\leq cr^s \textrm{ for all } x\in X, r>0&\}.
\end{align*}

$\m^s(X)\neq\emptyset \Rightarrow \h^s(X)>0$, indeed, if $\mu\in\m^s(X)$ and $E_1,E_2,\ppp$ is any covering of $X$, picking $x_j\in E_j, j=1,2,\ppp$ we have $B(x_j,d(E_j))\supset E_j$ and therefore
$$\sum_jd(E_j)^s\geq\sum_j\inv{c}\mu(B(x_j,d(E_j)))\geq\sum_j\inv{c}\mu(E_j),$$
 hence $\h^s(X)\geq\mu(X)/c>0$. The converse is also true for analytic sets in complete separable metric space, see for instance Theorem 8.17 in \cite{Mattila_1995}. We give a elementary proof for the same result starting from Howroyd \cite{Howroyd_1995}.
\begin{lem}\label{lem_Frostman}
If $X$ is an analytic subset of a complete separable metric space, then $\h^s(X)>0\Rightarrow\m^s(X)\neq\emptyset$.
\end{lem}
\begin{proof}
By Corollary 7 in \cite{Howroyd_1995}, we may assume $X$ compact with $0<\h^s(X)<\infty$. Let $M>10^s$ and $\delta>0$ such that $\h_\delta^s(X)>10^s\h^s(X)/M$ and consider the family 
$$\ca{B}:=\set{B(x,r): 10r\leq\delta, \h^s(B(x,r))\geq Mr^s}.$$ By the Vitaly covering theorem we get disjoint balls $B(x_j,r_j)\in\ca{B}, j=1,2,\ppp$ such that $\cup_{B\in\ca{B}}B\subset \cup_j B(x_j,5r_j)=:A$, hence
$$\h^s_\delta(A)\leq10^s\sum_j r_j^s\leq\frac{10^s}{M}\sum_j\h^s(B(x_j,r_j))\leq10^s\h^s(X)/M$$
and therefore $\h^s(X\backslash A)\geq\h^s_\delta(X\backslash A)\geq\h_\delta^s(X)-\h_\delta^s(A)>0$. Obviously, $\h^s$ restricted to $X\backslash A$ is in $\ca{M}^s(X)$ for $c=\max(M, 10^s\h^s(X)/\delta^s)$. 
\end{proof}

Let $\mu$ be a finite Borel measure on $X$. The \emph{$s$-energy} of $\mu$ is $$I_s(\mu)=\int\int\frac{d\mu(x_1)d\mu(x_2)}{d(x_1,x_2)^s}.$$

Finiteness of energy is closely relate with measure in $\m^s(X)$, for instance, if $I_s(\mu)<\infty$, then $\sup_{r>0}\abs{\mu(B(x,r))}/r^s\leq\int d(x,x_2)^{-s}d\mu(x_2)<\infty$ for $\mu$-a.e. $x\in X$. If $M$ large enough such that the Borel set $B=\set{x: \sup_{r>0}\abs{\mu(B(x,r))}/r^s\leq M}$ has positive $\mu$-measure, then $\nu=\mu|B$ is in $\m^s(X)$, in fact if $B(x,r)\cap B=\emptyset$ then $\nu(B(x,r))=0$, otherwise, if $z\in B(x,r)\cap B$, then $B(x,r)\subset B(z,2r)$ and therefore
$$\nu(B(x,r))\leq\nu(B(z,2r))\leq 2^sMr^s.$$

 On the other hand if $\mu\in\m^s(X)$, then $I_t(\mu)<\infty$ for all $0\leq t<s$, in fact
\begin{align*}
&\int d(x_1,x_2)^{-t}d\mu(x_2)=\int_0^\infty\mu[x_2:d(x_1,x_2)^{-t}\geq u]du\\
=&\int_0^{\infty}\mu(B(x_1,u^{-\inv{t}}))du\leq \mu(X)+c\int_1^\infty u^{-\frac{s}{t}}du<\infty.
\end{align*}

Then, if $X$ is an analytic subset of a complete separable metric space, we have:
\begin{align}\label{eq_dimension}
\dim X:=&\sup\set{s: \h^s(X)>0}\\\nonumber=&\sup\set{s:\exists\mu \textrm{ Borel measure with } 0<\mu(X)<\infty \text{ and }I_s(\mu)<\infty}.
\end{align}

\subsection{A scheme to compare measures.} We present here a elementary alternative to the methods of differentiation theory of Radon measures in Euclidean spaces that actually can not be extended to the general setting we are interested.  

Given any Borel measure $\nu$ on $Y$ we define the \emph{$\kappa$-density} as
$$D^{\kappa}\nu(y):=\liminf_{r\to0}\nu(B(y,r))/r^\kappa.$$

\begin{lem}\label{lem_key}
Let $\nu_1$ and $\nu_2$ be Borel measures, $\nu_2$ upper regular and $0<s,t<\infty$. If $D^{\kappa}\nu_1(y)\leq s$ and $D^{\kappa}\nu_2(y)\geq t$ for all $y\in A$, then $\nu_1(A)\leq 5^{\kappa}st^{-1}\nu_2(A).$ 
\end{lem}
\begin{proof}
For any $\e>0$ and any $U\supset A$ open set, by the Vitali covering theorem we get disjoint balls $B(x_j,r_j)\subset U$ with $\nu_1(B(x_j,5r_j))\leq (s+\e)(5r_j)^\kappa$ and $\nu_2(B(x_j,r_j))\geq (t+\e)r_j^\kappa$ for $j=1,2,\ppp$ such that $A\subset \cup_j B(x_j,5r_j)$, then
$$\nu_1(A)\leq\sum_j\nu_1(B(x_j,5r_j))\leq5^\kappa \frac{s+\e}{t-\e}\sum_j\nu_2(B(x_j,r_j))\leq 5^\kappa \frac{s+\e}{t-\e}\nu_2(U).$$
The inequality now follows from the upper regularity of $\nu_2$ and let $\e\to0$.
\end{proof}

This lemma has a strong version for Radon measure in Euclidean spaces, see \cite{Mattila_1995} Lemma 2.13, where for instant the $5^\kappa$-factor does not appear. 

\begin{pro}\label{pro_new}
Let $\nu$ and $\eta$ be Borel measures, $\eta$ upper regular and $D^{\kappa}\nu(y)<\infty$ for $\nu$-a.e. $y\in Y$.

\begin{itemize}
\item[(i)] If $D^{\kappa}\eta(y)>0$ for all $y\in Y$, then $\nu\ll\eta$.
\item[(ii)] If $D^{\kappa}\eta(y)\geq c>0$ for all $y\in Y$ and $\eta$ is $\sigma$-finite, then $d\nu/d\eta\leq 5^{\kappa}c^{-1}D^{\kappa}\nu$ for $\eta$-a.a.
\end{itemize} 
\end{pro}
\begin{proof}
(i) Let $A\subset Y$ with $\eta(A)=0$. Let $P_M=\set{y: D^{\kappa}\nu(y)<M,  D^{\kappa}\eta(y)>M^{-1}}.$ By Lemma \ref{lem_key} $\nu(A\cap P_M)\leq 5^{\kappa}M^2\eta(A\cap P_M)=0$. Let $M\to\infty$ we get $\nu(A)=0$. 

(ii)  Let  $B\subset Y$ any Borel set. For $1<t<\infty$ let $B_p=\set{y\in B: t^p\leq D^{\kappa}\nu(y)< t^{p+1}}$, $p\in\Z$. By the Lemma \ref{lem_key} $\nu(B_p)\leq 5^\kappa t^{p+1} c^{-1}\eta(B_p)$ and $\nu(\set{D^{\kappa}\nu=0})$, then
\begin{align*}
\nu(B)&=\sum_{p\in\Z}\nu(B_p)\leq \sum_{p\in\Z}5^\kappa t^{p+1} c^{-1}\eta(B_p)\\&\leq t5^\kappa c^{-1} \sum_p\int_{B_p}D^{\kappa}\nu d\eta\leq t5^\kappa c^{-1} \int_{B}D^{\kappa}\nu d\eta.
\end{align*}
Let $t\to 1$, we get $\nu(B)\leq 5^{\kappa}c^{-1}\int_BD^{\kappa}\nu d\eta$ for all Borel set $B\subset Y$. 
\end{proof}

This Proposition is in some sense related with Theorem 2.12 in \cite{Mattila_1995}.

\section{Proof of the Theorems}\label{section_results}
Let $X,Y$ be separable metric spaces, $(\Lambda, \p)$ a probability space and $\pi:\Lambda\times X\to Y$ a measurable function satisfying the transversality property in \ref{transversality}. 

Given any finite Borel measure $\mu$ on $X$ let $\nu_\la=(\pi_{\la})_*\mu$. Note that $\nu_\la$ are also Borel measures with $\nu_\la(\pi_\la(X))=\mu(X)$.

The following two lemmas are coming naturally from ideas of Theorem 9.3 and Theorem 9.7 in \cite{Mattila_1995}.

\begin{lem}\label{thm1}
$$\int_{\Lambda}I_t(\nu_\la)d\p\leq C_{t,\kappa}I_{\alpha t}(\mu) \text{ for all } 0\leq t<\kappa.$$
\end{lem} 
\begin{proof} Notice that $I_t(\nu_\la)=\int\int\frac{d\mu(x_1)d\mu(x_2)}{d(\pi_\la(x_1),\pi_\la(x_2))^t}$ then, by Fubuni's theorem we have
$$\int_{\Lambda}I_t(\nu_\la)d\p=\int\int\Big[\int_{\Lambda}\paren{\frac{d(\pi_\la(x_1),\pi_\la(x_2))}{d(x_1,x_2)^\alpha}}^{-t} d\p \Big]\frac{d\mu(x_1)d\mu(x_2)}{d(x_1,x_2)^{\alpha t}},$$
and 
\begin{align*}
\int_{\Lambda}\paren{\frac{d(\pi_\la(x_1),\pi_\la(x_2))}{d(x_1,x_2)^\alpha}}^{-t} d\p&=\int_{0}^{\infty}\p\Big[\la\in\Lambda: \paren{\frac{d(\pi_\la(x_1),\pi_\la(x_2))}{d(x_1,x_2)^\alpha}}^{-t}\geq u\Big]du\\
&=\int_{0}^{\infty}\p\Big[\la\in\Lambda: \frac{d(\pi_\la(x_1),\pi_\la(x_2))}{d(x_1,x_2)^\alpha}\leq u^{-\inv{t}}\Big]du\\
&\leq 1+ C\int_{1}^{\infty}u^{-\frac{\kappa}{t}}du<\infty.
\end{align*}
\end{proof}

\begin{proof}[Proof of Theorem \ref{thm_Marstrand1}]
It follows from equation (\ref{eq_dimension}) and Lemma \ref{thm1}. 
\end{proof}

\begin{lem}\label{thm2}
$$\int_\Lambda\int D^{\kappa}\nu_\la d\nu_\la d\p\leq C I_{\alpha\kappa}(\mu).$$
\end{lem}
\begin{proof} Using Fatou's lemma 
$$\int\liminf_{r\to0}\frac{\nu_\la(B(y,r))}{r^\kappa}d\nu_\la(y)\leq \liminf_{r\to0}\inv{r^\kappa}\int\nu_\la(B(y,r))d\nu_\la(y), \text{ and}$$
$$\int\nu_\la(B(y,r))d\nu_\la(y)=\nu_\la\times\nu_\la[(y_1,y_2): d(y_1,y_2)\leq r]=\mu\times\mu[(x_1,x_2): d(\pi_\la(x_1),\pi_\la(x_2))\leq r],$$
then, by Fubini's theorem and (\ref{transversality})
\begin{align*}
\int_\Lambda\int D^{\kappa}\nu_\la d\nu_\la d\p
&\leq\liminf_{r\to0}\iint\inv{r^\kappa}\p[\la\in \Lambda: d(\pi_{\la}(x_1),\pi_{\la}(x_2))\leq r]d\mu(x_1)d\mu(x_2)\\&\leq C I_{\alpha\kappa}(\mu).
\end{align*}
\end{proof}

\begin{proof}[Proof of Theorem \ref{thm_Marstrand2}]
(i) follows from Lemma \ref{thm2}, Lemma \ref{lem_Frostman} and Proposition \ref{pro_new}(i), since if $\nu_\la\ll\eta$ with $\nu_\la(\pi_\la(X))=\mu(X)>0$ then $\eta(\pi_\la(X))>0$. 

Analogously with Proposition \ref{pro_new}(ii) we have   $d\nu_\la/d\eta\in L^2(\eta)$ for a.e. $\la\in\Lambda$ with 
$$\int_\Lambda\norm{\frac{\nu_\la}{d\eta}}^2_{L^2({\eta})}d\p\leq 5^\kappa c^{-1}C I_{\alpha\kappa}(\mu)<\infty.$$
The part (ii) now follows from the Cauchy-Schwarz's inequality $\mu(X)\leq\eta(\pi_\la(X))\norm{\frac{\nu_\la}{d\eta}}^2_{L^2({\eta})}$.
\end{proof}

\bibliographystyle{amsplain}
\bibliography{bibli}

\providecommand{\bysame}{\leavevmode\hbox to3em{\hrulefill}\thinspace}
\providecommand{\MR}{\relax\ifhmode\unskip\space\fi MR }
\providecommand{\MRhref}[2]{%
  \href{http://www.ams.org/mathscinet-getitem?mr=#1}{#2}
}
\providecommand{\href}[2]{#2}
\begin{thebibliography}{1}

\bibitem{Howroyd_1995}
J.D. Howroyd, \emph{On dimension and on the existence of sets of finite
  positive {Hausdorff} measure}, Proc. London Math. Soc. (3) \textbf{70}
  (1995), 581--604.

\bibitem{Mattila_1995}
P.~Mattila, \emph{Geometry of sets and measures in euclidean spaces : fractals
  and rectifiability}, Cambridge University Press, 1995.

\bibitem{Solomyak_1998}
Boris Solomyak, \emph{Measure and dimension for some fractal families}, Math.
  Proc. Cambridge Philos. Soc. \textbf{124} (1998), no.~3, 531--546.

\end{thebibliography}
\end{document}